\documentclass[11pt]{article}
\usepackage{amsmath}
\usepackage{amsfonts}
\usepackage{authblk}
\usepackage{url}
\usepackage[english]{babel}
\usepackage[dvipsnames]{xcolor}
\usepackage{hyperref}
\usepackage{graphicx}
\usepackage{caption}
\usepackage{subcaption}
\usepackage[utf8]{inputenc}
\usepackage{booktabs}
\usepackage{changes}
\usepackage{latexsym,amsthm,amscd}
\newcounter{alphthm}
\usepackage{tikz}
\usetikzlibrary{shapes.geometric, arrows}
\tikzstyle{startstop} = [rectangle, rounded corners, minimum width=1cm, minimum height=1cm,text centered, draw=black, fill=blue!30]
\tikzstyle{io} = [trapezium, trapezium left angle=70, trapezium right angle=110, minimum width=3cm, minimum height=1cm, text centered, draw=black, fill=blue!30]
\tikzstyle{process} = [rectangle, minimum width=5cm, minimum height=1cm, text centered, draw=black, fill=white!30]
\tikzstyle{decision} = [ellipse, minimum width=3cm, minimum height=1cm, text centered, draw=black, fill=white!30]
\tikzstyle{arrow} = [thick,->,>=stealth]
\tikzstyle{edge from parent}=[draw,dashed,thick,->,>=stealth,black]
\tikzstyle{process} = [rectangle, minimum width=3cm, minimum height=1cm, text centered, text width=3cm, draw=black, fill=white!30]
\theoremstyle{remark}
\usepackage{amssymb, graphicx, tcolorbox, wasysym,  pifont}
\numberwithin{equation}{section}

\newcommand{\rev}[1]{\textcolor{black}{#1}}

\newtheorem*{theorem*}{Theorem} 
\newtheorem{theorem}{Theorem}[section]

\theoremstyle{definition}

\newtheorem{example}[theorem]{Example}

\theoremstyle{remark}
\newtheorem{remark}[theorem]{Remark}
\renewcommand{\theequation}
{\arabic{section}
.\arabic{equation}}
\title{A Geometric Framework for Wavefront-Based Wildfire Spread Modeling on  Inclined Landscapes}

\date{\today}

\author{Hengameh R. Dehkordi}
\affil{Center of Mathematics, Computing, and Cognition, Federal University of ABC, 09210-580, Santo André, SP, Brazil\\
\texttt{hengameh.r@ufabc.edu.br}}

\begin{document}
\maketitle

\begin{abstract}
\rev{This paper develops a geometric framework for the propagation of wildfire fronts on locally inclined terrain. The framework connects prescribed head-fire and backing-fire rates to an explicit directional wavelet, with wind and uphill slope combined through a vector-weighting rule. A common construction covers homogeneous and heterogeneous media, both with and without wind. The resulting wavelet reproduces the specified forward and backward rates and admits windless slope-driven and flat-terrain limits. We establish the positivity and strong-convexity requirements for a regular Finsler indicatrix, including an exact threshold for the windless cosine profile and an angular criterion for the general wind--slope case. Front evolution is characterized by minimizing outward Finsler-normal geodesics up to their cut times, or equivalently by the boundary of a forward-distance reached region. A piecewise-static approximation accommodates updates to environmental inputs. Numerical illustrations with nonconvex initial fronts, discretization refinement, and multi-phase grid propagation demonstrate the geometric effects of anisotropic local wavelets. The prescribed rates and other numerical inputs are synthetic; their empirical estimation, field validation, and comparison with operational simulators remain subjects for future applied studies.}

{\textcolor{white}{Mathematical Methods in the Applied Sciences}}
\end{abstract}


\section{Introduction}

Wildfires pose a major threat to forests, grasslands, wildlife, and human infrastructure, often causing severe ecological and economic damage. Understanding and predicting fire behavior, particularly the evolution of the wildfire perimeter, is therefore essential for effective wildfire management and hazard mitigation. Predicting perimeter movement remains challenging because fire spread is a dynamic process influenced by fuel characteristics, landscape topography, and wind \cite{andrews2018rothermel,speer2022wildland}. Wind can tilt flames toward unburned vegetation and increase preheating, while wind and terrain conditions vary across space and time \cite{sharples2009overview,dai1999diurnal}.

Wildfire spread models must account for the physical processes and environmental factors that influence fire behavior. However, not all these processes can be described explicitly at operational scales. One of the main challenges in wildfire modeling is to establish a mathematically consistent link between local fire spread rates and the evolution of large-scale fire perimeters. For example, the widely used Rothermel surface fire model estimates the local rate of spread under specified fuel, moisture, wind, and slope conditions. This estimate is commonly used to determine the maximum rate of spread in the direction of greatest fire advance. However, the model does not, by itself, describe how the spread rate varies with direction or how a two-dimensional fire front evolves. These aspects require additional assumptions and a suitable front propagation framework \cite{rothermel1972mathematical,andrews2018rothermel}.

Several mathematical models and numerical simulators have been developed to describe fire front propagation. However, an important challenge is to derive an explicit direction dependent rate of spread from prescribed head-fire and backing-fire rates. This requires a mathematical formulation that describes how the spread rate varies with direction while satisfying the necessary mathematical conditions. The formulation should also account for the combined effects of wind and terrain slope, including cases where they act in different directions. It should be applicable to both homogeneous and spatially heterogeneous landscapes.

To narrow this gap, this work develops a simplified  geometric framework of fire front propagation on locally sloped terrain. In particular, we determine a head-fire rate and a backing-fire rate as input to the model and construct an explicit directional wavelet. By the head-fire rate, \(R_H\), we mean the rate at which the fire front advances in the direction of maximum spread. The backing-fire rate, \(R_B\), is the rate of advance in the opposite direction. Terrain effects are represented by the slope magnitude and the uphill direction.  This wavelet exactly recovers these rates in the preferred forward and opposite backward directions. The directional effects of wind and uphill slope are incorporated by presenting a wind-slope weighting rule. This directional wavelet is employed in a  geometric propagation and a Huygens-type region construction to model global front evolution. Rather than replacing detailed combustion models or operational simulators, the proposed approach complements them by providing an explicit, mathematically verifiable rate-to-wavelet formulation together with its governing admissibility conditions.
Our approach builds on two main areas of research: practical fire spread models and geometric theories of wavefront propagation.

Richards' seminal equations in the operational and empirical literature describe the evolution of fire fronts using expanding elliptical wavelets coupled with a system of differential equations \cite{richards1990elliptical}. Operational simulators such as Phoenix, IGNITE, Bushfire, FireMaster, FARSITE and Prometheus employ a variety of fire behaviour and perimeter propagation procedures \cite{Sullivan2009}. Most notably, FARSITE integrates local fire behavior estimates with Huygens type wavelet  to propagate the fire perimeter \cite{farsite2004fire}. The elliptical wavelet is an additional directional modeling assumption, and is not a direct consequence of the Rothermel spread rate equation. Operational frameworks may also allow for dynamic updates of weather and environmental inputs during a simulation, meaning that local wavelet propagation does not need to assume fixed environmental conditions for the entire duration of a wildfire event.

Another paradigm for front tracking is given by implicit representations, for example level-set methods. Level-set methods allow the representation of a moving front as an implicit function, as well as natural topological changes and merging of fronts. Wavelet based Huygens constructions give an explicit description of the local growth by means of directional wavelets and the boundary of the reached region. Both paradigms need a prescribed law of direction of propagation. It does not automatically set such things as fuel, weather or combustion parameters. In the meantime, Finsler geometry has been an important mathematical tool for modeling anisotropic wavefronts and their optimal propagation paths. In \cite{markvorsen2016finsler} the application of Huygens' principle was studied in the case of two-dimensional Finsler spaces with elliptical frames. This basis was later generalized to a geodesic-spray interpretation of fire front propagation, which made the connection to Richards' equations and time-dependent wind conditions clear \cite{markvorsen2017geodesic}. For extreme propagation regimes with strong winds, where the upwind spread may vanish, geometric formulations based on Lorentz-Finsler geometry with cones have been introduced \cite{javaloyes2021applications} and related geometric developments for sloped terrain have been reported \cite{angel2023general}. Further theoretical developments have led to wave equations for wind-driven fire spread on flat terrain \cite{dehkordi2022applications} and to the study of strategic paths for wildfire management and infrastructure protection \cite{dehkordi2023applications,dehkordi2023wildfire,dehkordi2021mathematical} that are closely related to the larger mathematical blocking problem \cite{bressan2013}.

While these contributions establish essential theoretical foundations, the present study builds upon and complements them by developing an explicit directional wavelet model that directly incorporates prescribed head-fire and backing-fire rates, together with wind and terrain characteristics, to describe wildfire-front propagation under different environmental conditions. The proposed framework provides a practical mathematical formulation that can be applied to real-world scenarios using available environmental data, while explicitly characterizing the admissibility conditions under which the resulting directional wavelet defines a regular Finsler metric. In this respect, the main contributions of this paper can be summarized as follows:

\begin{enumerate}
\item {Construction of rates to wavelets explicitly:} We derive an analytical directional propagation law that is directly parameterized from given head-fire ($R_H$) and backing-fire ($R_B$) rates and equipped with an explicit vector weighting of wind and uphill directions. The proposed construction recovers the prescribed forward and backward rates, and reduces to windless and flat-terrain limits. 
\item {Mathematical admissibility conditions:} We study the necessary conditions of positivity and strong convexity on the local directional wavelet to define a regular geometric structure. In particular, we distinguish the exact windless cosine-wavelet admissibility threshold from the windy case where a full angular convexity condition must hold. Analytical derivations for special cases are augmented with systematic numerical angular checks for the combined wind-slope construction. \item {A single propagation framework for four environmental settings:} In one geometric model we explore four different model configurations: homogeneous and spatially heterogeneous media, with and without wind. We show the effect of the local wavelet geometry on the evolving fire perimeter by means of the Huygens-type front evolution and with demonstrative numerical experiments. We study the sensitivity to discretisation by numerical grid refinement.
\end{enumerate}

It is useful to note some intentional modeling limitations in defining the scope. We do not include phenomena such as firebrand spotting and active suppression. Rather than solving coupled partial differential equations for radiation, convection, conduction or atmospheric feedback, we include the directional effects of wind and slope using a local kinematic wavefront law in our framework. In extreme wind conditions, when the upwind spread is totally suppressed and the directional speed is null, the propagation is outside the positive, strongly convex Finsler setting discussed here, and a cone-based approach would be required. Such phenomenon is extremely rare in real-world. Also, for mathematical tractability on locally plane terrain, the environmental coefficients are time-invariant within each propagation interval, but they can be space-varying and updated between intervals. This piecewise-static approximation does not imply that the real wildfire environments are stationary; its temporal approximation error is not quantified here. Finally, because this work focuses on the mathematical rate-to-wavelet closure and its geometric properties, we treat $R_H$ and $R_B$ as prescribed inputs; parameter estimation from remote-sensing data, empirical calibration of the ratio $R_B=kR_H$, $0<k<1$, field validation, and benchmarking against level-set or operational simulators are left for future applied studies.

The rest of the paper is organized as follows. In Section 2 we recall the necessary preliminaries in Finsler geometry. Section 3 develops the rate-to-wavelet construction for homogeneous and heterogeneous media, with and without wind. Section 3.2 discusses geodesic cut times and numerical limitations, and Section 3.3 summarizes the algorithmic flow of the proposed construction. Some geometric illustrations are given in Section \ref{examp} as well as the numerical implementation, a discretization refinement analysis and an example of propagation on a grid. We conclude the paper in Section \ref{concl}. Finally, in Appendix \ref{apen3} we describe the rate formulation and the recovery of the head-fire and backing-fire rates. Appendix \ref{geodesics} discusses the mathematical conditions for admissibility and their numerical evaluation.

\section{Preliminaries}

We consider a locally inclined terrain surface $M:z=h(x_1,x_2)$ approximated locally by its tangent plane, projected onto horizontal coordinates
\(\{(x_1,x_2,z)\in\mathbb R^3:z=d_1x_1+d_2x_2+d\}\). A point of \(M\)
is therefore represented by \((x_1,x_2,d_1x_1+d_2x_2+d)\), and a tangent
vector is identified with its horizontal components \(v=(v^1,v^2)\). The metric induced on \(M\) by the ambient Euclidean space is
\[
ds^2=(1+d_1^2)\,dx_1^2+2d_1d_2\,dx_1\,dx_2+(1+d_2^2)\,dx_2^2.
\]
 Moreover, at each point $(x_1,x_2)$, the slope is $S = \sqrt{d_1^2 + d_2^2} = \tan\phi$ and the uphill direction is $\Theta_s = \operatorname{atan2}(d_2,d_1)$. Throughout this paper
\(\|v\|=\sqrt{(v^1)^2+(v^2)^2}\) denotes the Euclidean norm of the
{horizontal projection}.
A Finsler metric is a function \(F:TM\to[0,\infty)\) that is smooth on
the tangent bundle \(TM\setminus\{0\}\), positive for nonzero
tangent vectors, and positively homogeneous of degree one:
\(F(x,\lambda v)=\lambda F(x,v)\) for every \(\lambda>0\).
It is also required to be strongly convex, meaning that its fundamental
tensor
\begin{equation}\label{sff}
 g_{ij}(x,v)=\frac12\frac{\partial^2F^2}{\partial v^i\partial v^j}(x,v),
 \qquad i,j=1,2,
\end{equation}
is positive definite for every \(v\ne0\). 

The  geodesics are the shortest time-minimizing paths and they satisfy
\begin{equation}\label{geo-fin}
\frac{d^2x^i}{dt^2}+2G^i\left(x,\frac{dx}{dt}\right)=0,
\qquad i=1,2,
\end{equation}
where 
\[
G^i(x,v)=\frac14g^{i\ell}(x,v)
\left(\frac{\partial^2F^2}{\partial x^k\partial v^\ell}(x,v)v^k
-\frac{\partial F^2}{\partial x^\ell}(x,v)\right).
\]
Here \(g^{ij}\) denotes the inverse of the fundamental tensor
\(g_{ij}\).

A wavelet represents the boundary of the region reached by a fire spreading from a single ignition point over a specified time interval.
The burned region at time \(t\) is denoted by \(K_t\), and its fire
front is \(\partial K_t\). A fire ray is a propagation path associated
with a geodesic of the Finsler metric. When the metric is well defined,
a point ignition produces the wavelet \cite{dehkordi2022applications, dehkordi2018finsler, alexandrino2019finsler}. The set of Finsler unit  vectors at a
point is called the indicatrix.  Under locally homogeneous conditions, the wavelet for unit propagation time coincides with the indicatrix. 

The following theorem provides a geometric basis for describing fire
fronts by Finsler distance level sets. Here \(K_0\) is a compact initial
set, and the stated no-cut-locus assumption applies to the region
through which the fronts evolve.

\begin{theorem}\label{one} \cite{dehkordi2019huygens}
Let $(M,F)$ be a Finsler manifold on which a wildfire propagates, and let $K_0 \subset M$ be compact. Define the distance function $\rho:M\to\mathbb{R}$ by
\[
\rho(x)=d_F(K_0,x), 
\] where  $d_F(K_0,x)$ is the forward Finsler distance from $K_0$ to $x$. 
Suppose that $\rho^{-1}([s,r])=\hat{U}$ and that $\rho^{-1}((s,r))$ contains no cut points of $K_0$. If $\rho^{-1}(s)$ represents the fire front at time $0$, then $\rho^{-1}(t)$ represents the fire front at time $t-s$, for every $t\in[s,r]$. Moreover, the family of fronts
\[
\{\rho^{-1}(t)\}_{t\in[s,r]}
\]
satisfies Huygens' principle. Each fire ray is a geodesic that is orthogonal, in the Finsler sense, to the fire fronts it intersects.
\end{theorem}

In this theorem, orthogonality is understood in the Finsler sense,
through the fundamental tensor, rather than necessarily as Euclidean
orthogonality in the horizontal plane. The no-cut-locus assumption
ensures that the relevant geodesic segments remain minimizing and the
distance level sets provide regular advancing fronts. More generally,
the burned region is described by the forward-distance sublevel set;
this reached region formulation remains appropriate when distinct rays
meet and the front develops nonsmooth points.

Huygens' principle describes the evolution of a wavefront through
secondary point-source wavefronts. In the terminology of forward
Finsler distance, the classical statement can be expressed as follows.

\begin{theorem}\cite{arnol2013mathematical}
Let $\phi_x(t)$ be the spherical wavefront of the point $x$ after time $t$. For every point $y\in\phi_x(t)$, consider the spherical wavefront after time $s$, i.e. $\phi_y(s)$. Then, the spherical wavefront of point $x$ after $s+t$ will be the envelope of spherical wavefronts $\phi_y(s)$.
\end{theorem}

For fire fronts, the envelope is understood as the boundary of the
reached region, rather than the unrestricted collection of propagated
ray endpoints. The next section introduces the directional wavelets,
states the positivity and strong convexity conditions under which they
define regular Finsler metrics, and develops the corresponding front
and ray constructions for the model configurations.

\section{Propagation on an Locally Inclined Terrain}
This section develops the geometric description of wildfire propagation on an locally inclined terrain under spatially homogeneous and heterogeneous environmental conditions, with or without wind. All propagation velocities and distances are expressed in the horizontal projected coordinates introduced in the Preliminaries. Throughout, \(\langle\cdot,\cdot\rangle\) denotes the standard Euclidean inner product in these coordinates.

\paragraph{Common geometric notation.}

Write \(e_\theta=(\cos\theta,\sin\theta)\) and \(e_\theta^\perp=(-\sin\theta,\cos\theta)\). For a positive directional propagation speed $r(x,\theta)$,
 define the
filled local unit-time wavelet at $x$ by
\begin{equation}\label{waveletbody}
\overline{\mathcal{I}_x}
=
\{\rho e_\theta:
0\leq\rho\leq r(x,\theta),\
0\leq\theta<2\pi\}.
\end{equation}
Its boundary is the indicatrix
\[
\mathcal{I}_x
=
\{r(x,\theta)e_\theta:
0\leq\theta<2\pi\}.
\]
This construction represents the local propagation geometry
with the environmental conditions fixed at $x$.
For spatially heterogeneous conditions, the evolving fire
front is determined by the corresponding position-dependent
propagation model.
Thus, \(\overline{\mathcal{I}_x}\) is the filled set of locally admissible unit-time displacement vectors, whereas \(\mathcal{I}_x\) is the corresponding local unit wavefront. The translated and scaled set \(x+\Delta t\,\overline{\mathcal{I}_x}\) is the exact point-source reached region over \(\Delta t\) when the metric is spatially homogeneous and time independent on that interval. In a spatially heterogeneous medium, freezing the coefficients at \(x\) yields only a local approximation to the forward-distance wavelet.

For the directional law to define a regular, strongly convex  indicatrix, we require, at every relevant position \(x\) and direction \(\theta\), see condition (A2) in Appendix~\ref{geodesics},
\[
r(x,\theta)>0,\qquad
r(x,\theta)^2+2[\partial_\theta r(x,\theta)]^2
-r(x,\theta)\partial_{\theta\theta}r(x,\theta)>0.
\]
Here \(r\) is assumed sufficiently smooth in its angular variable. These conditions ensure that the indicatrix is a smooth, positively curved boundary enclosing the origin.

At a smooth point \(x\in\partial K_t\), let \(n_x\) denote the outward Euclidean unit normal. The outward Finsler-unit normal velocity is the supporting point of the filled indicatrix in direction \(n_x\):
\begin{equation}\label{supportnormal}
v_x=\operatorname*{arg\,max}_{w\in\overline{\mathcal{I}_x}}(n_x\cdot w),
\end{equation}
where \(\overline{\mathcal{I}_x}=\{w:F(x,w)\leq1\}\).
\begin{figure}[h!]
			\centering
                   \includegraphics[width=\textwidth]{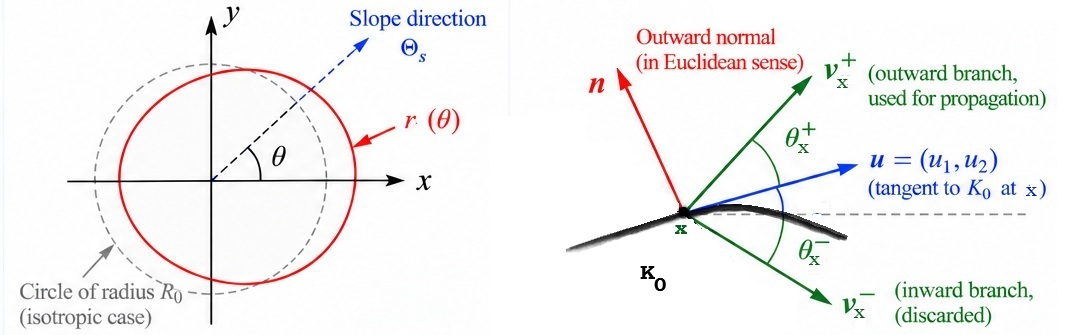}
         \caption{The normal vectors and the indicatrix.}\label{waves}
         \end{figure}
Here, \(v_x\) denotes the outward propagation velocity vector, whose direction and magnitude are determined by the supporting point of the local wavelet in the direction of the outward normal covector \(n_x\), see Fig.\ref{waves}. Positivity and strong convexity ensure that this maximizer is unique and lies on \(\mathcal{I}_x\).
This supporting-point construction yields Finsler orthogonality: since \(F(x,v_x)=1\), the Lagrange multiplier condition gives
\[
n_x=\eta\,d_vF(x,v_x),\qquad \eta>0,
\]
where $\eta$ is a scalar ensuring that the two covectors are positively proportional. The Euclidean inner product identifies the normal vector with a covector. For every tangent vector \(u\in T_x\partial K_t\), we have \(n_x\cdot u=0\), and therefore \(d_vF(x,v_x)[u]=0\). Using the homogeneity identity
\[
g_{x,v_x}(v_x,u)=F(x,v_x)\,d_vF(x,v_x)[u],
\]
we obtain
\[
g_{x,v_x}(v_x,u)=0.
\]
Thus, the supporting-point construction identifies the unique outward Finsler-unit normal velocity without directly solving the orthogonality equation.

The radial propagation speed \(r(x,\theta)\) must be distinguished from the normal advance speed of the front. In an anisotropic medium, \(v_x\) need not be parallel to \(n_x\).  The Euclidean normal specifies the supporting direction, whereas \(v_x\) gives the actual Finsler-normal propagation velocity.

\paragraph{Piecewise-static temporal setting.}
Let \(0=t_0<t_1<\cdots<t_N=T\) be prescribed environmental update times. During each interval \(I_k=[t_k,t_{k+1})\), the environmental parameters are fixed in time and determine the Finsler metric
\[
F_k(x,v)=\frac{\lVert v\rVert}{ r(x,\theta)},\qquad v\ne0,\quad \theta=\arg(v),
\]
where the radial speed function \( r(x,\theta)\) satisfies the strong convexity condition stated below. The metric may depend on position \(x\), but not on time within \(I_k\). When the forward distance construction is applicable, the burned region at time \(t_{k+1}\) is
\[
K_{t_{k+1}}=\{y\in\mathbb R^2:d_{F_k}(K_{t_k},y)\leq\Delta t_k\},\qquad
 d_{F_k}(K_{t_k},y)=\inf_{x\in K_{t_k}}d_{F_k}(x,y),
\] where  \(\Delta t_k=t_{k+1}-t_k\).
The reached region at the end of an interval becomes the initial region for the next. This defines the piecewise-static mathematical model; the finite-step unions of sampled local wavelets used in the numerical implementation (Section~\ref{algsec}) approximate these distance-based updates. The construction represents temporal changes in environmental conditions through a piecewise-static approximation. If all the metrics \(F_k\) coincide, it reduces to the time-independent setting.

The head-fire and backing-fire rates must satisfy \(R_H\geq R_B>0\) with \(R_B=kR_H\) for a prescribed ratio \(0<k\leq1\). These rates, together with wind and terrain information, enter the phenomenological directional-wavelet construction described below and in Appendix~\ref{apen3}. The paper does not derive an empirical law for \(k\), nor does the geometric formulation resolve the underlying combustion and heat-transfer processes, as these are beyond the slope of this work.

\subsection{Wildfire propagation model}
Throughout this section, $R_H$ and $ R_B$  are prescribed inputs 
which may be inferred from suitable observations, see for instance \cite{dehkordi2026satellite, dehkordi2026determining}, moreover we have $R_H \geq R_B>0$.

\subsubsection{Modeling general wildfire propagation in homogeneous or heterogeneous media, with or without wind.}
Suppose $0<R_B(x)\le R_H(x)$, wind speed $U(x)\ge0$,
wind direction $\Theta_w(x)$ wherever $U(x)>0$, and terrain slope angle $\phi(x)\ge0$ are prescribed within a time-frozen interval. Define the
dimensionless wind speed $U^*(x)=U(x)/U_{\rm ref}$, where
$U_{\rm ref}=1\,\mathrm{m/s}$, and set
\[
\begin{gathered}
b(x)=\frac{R_H(x)+R_B(x)}{2},\qquad
a(x)=\frac{b(x)}{1+0.25U^*(x)},\\
\lambda(x)=\frac{S(x)}{S(x)+\mu U^*(x)},\qquad
c(x)=(1-\lambda(x))\frac{R_H(x)-R_B(x)}{2},\\
g(x)=\lambda(x)\frac{R_H(x)-R_B(x)}{2},
\end{gathered}
\]
where $\mu>0$ is the prescribed dimensionless wind--slope weighting
parameter and $S(x)=\tan\phi(x)$ for terrain slope angle $\phi(x)$.
When $S(x)>0$ or $U(x)>0$, 
 the preferred propagation direction is
obtained from the wind--slope resultant
\[
\begin{gathered}
(q_x(x),q_y(x))=(1-\lambda(x))
(\cos\Theta_w(x),\sin\Theta_w(x))
+\lambda(x)(\cos\Theta_s(x),\sin\Theta_s(x)),\\
\Theta(x)=\operatorname{atan2}(q_y(x),q_x(x)),
\end{gathered}
\]
where $\Theta_s(x)$ is the uphill direction and $\Theta_w(x)$ is the
wind direction. At $S(x)=0$ with $U(x)>0$, use the flat-terrain limit
$\lambda(x)=0$ and $\Theta(x)=\Theta_w(x)$ without evaluating an uphill
angle. At $U(x)=0$ with $S(x)>0$, $\lambda(x)=1$ and
$\Theta(x)=\Theta_s(x)$ without evaluating a wind angle. If both wind
and slope vanish, the ratio defining $\lambda$ is undefined and a
separate zero-wind, flat-terrain formulation is required. We exclude a
vanishing wind--slope resultant, for which $\Theta$ is undefined.

Require the scalar coefficients $a,b,c,g$ and the unit preferred-direction
field $(\cos\Theta,\sin\Theta)$ to be $C^2$ in position on the modeled
region, with $a,b>0$ and $c<b$. These are sufficient spatial regularity
assumptions for the local Finsler spray, provided the angular admissibility
conditions below hold. Where $\Theta_s$ is obtained by normalizing a
terrain-elevation gradient, require that gradient to be nonzero; at a
zero-slope point the uphill angle is undefined and the formulation must use
its separately defined flat-terrain limit. The wind direction likewise
requires positive wind speed wherever it is used. Where the wind--slope
resultant is constructed from smooth input vector fields, assume
$\|(q_x(x),q_y(x))\|\ge\varepsilon>0$ on the region under consideration;
this excludes cancellation and ensures that normalizing the resultant
preserves $C^2$ regularity. If a spatial field crosses a zero-slope or
zero-wind set, the resulting coefficients and directional field must still
satisfy the stated $C^2$ conditions; piecewise smooth input data do not
automatically define a globally smooth Finsler spray.

The local directional propagation rate is
\begin{equation}\label{model4wavelet}
r(x,\theta)=\frac{c(x)a(x)^2\cos\delta_x+a(x)b(x)
\sqrt{a(x)^2\cos^2\delta_x+b(x)^2\sin^2\delta_x-c(x)^2\sin^2\delta_x}}
{a(x)^2\cos^2\delta_x+b(x)^2\sin^2\delta_x}
+g(x)\cos\delta_x,
\end{equation}
where $\delta_x=\theta-\Theta(x)$. 

\paragraph{Recovery of the prescribed rates.}
For the general radial law, setting $\delta_x=0$ or $\pi$ gives, respectively,
\[
r(x,\Theta(x))=b(x)+c(x)+g(x)=R_H(x),\qquad
r(x,\Theta(x)+\pi)=b(x)-c(x)-g(x)=R_B(x),
\]
because $c(x)+g(x)=[R_H(x)-R_B(x)]/2$. These are radial rates along the preferred direction and its opposite; they are not claims about the maximum front-normal velocity. The crosswind radial rate is
\[
r(x,\Theta(x)\pm\pi/2)=\frac{a(x)\sqrt{b(x)^2-c(x)^2}}{b(x)}.
\]
Matching the head and backing rates does not uniquely determine the directional profile. The offset-ellipse-plus-cosine form is a modeling choice that provides an explicit radial expression and controllable transverse width; it is not derived from the Rothermel model.

At each smooth point $x\in\partial K_0$ of the initial fire front,
let $n_x$ be the outward Euclidean unit normal. The outward
Finsler-unit propagation velocity is determined by
\begin{equation}\label{model3-support}
v_x=\operatorname*{arg\,max}_{w\in\mathcal I_x}n_x\cdot w.
\end{equation}
The corresponding fire ray $\gamma_x(t)$ is the geodesic with initial
data $\gamma_x(0)=x$ and $\dot\gamma_x(0)=v_x$. The burned region at time $\tau$ is
\begin{equation}\label{model3front}
\begin{aligned}
K_\tau=K_0\cup\bigl\{\gamma_x(t):\;&x\in\partial K_0,
\quad 0\le t\le\tau,\\
&\gamma_x|_{[0,t]}\text{ is minimizing from }K_0\bigr\}.
\end{aligned}
\end{equation}

\medskip\noindent\textit{Mathematical construction of Model 1.}
The position-dependent
Finsler metric is
\begin{equation}\label{fins4}
F_x(v)=\frac{\lVert v\rVert}{r(x,\arg v)}\quad(v\ne0),
\qquad F_x(0)=0.
\end{equation}
The filled indicatrix
$\overline{\mathcal I_x}=\{v:F_x(v)\le1\}$ is strongly convex under
the stated positivity and curvature conditions. Hence the maximizer
in Eq.~\eqref{model3-support} is unique at each smooth point of the
initial front. Initial data for the spatial geodesic equation
\eqref{geo-fin} satisfy the unit-speed and outward Finsler-normal
conditions
\begin{equation}\label{sph}
F_x(v_x)=1,
\end{equation}
\begin{equation}\label{sph2}
g_{x,v_x}(v_x,u)=0\qquad(u\in T_x\partial K_0).
\end{equation}
Indeed, $g_{x,v}(v,u)=F_x(v)\,d_vF_x(v)[u]$, and the supporting-point
condition makes $d_vF_x(v_x)$ proportional to the outward Euclidean
conormal. Writing $v_x=r(x,\theta)e_\theta$, these conditions become
\begin{equation}\label{model4normal}
\lVert v_x\rVert=r(x,\theta),\qquad
u\cdot[r(x,\theta)e_\theta-
\partial_\theta r(x,\theta)e_\theta^\perp]=0,
\end{equation}
where $u$ is any nonzero tangent vector to $\partial K_0$ at $x$.
The unique outward angular solution is fixed by the support construction.
Spatial derivatives enter the geodesic spray and do not enter this
pointwise angular condition; the spatial geodesic equations include all
necessary base-point derivatives. Between environmental updates the
metric is time independent; at an update the previous reached region
becomes the initial region for the next interval. No exactness claim
is made for finite-step frozen-wavelet propagation in heterogeneous media.

\medskip\noindent\textbf{Corollary 1 (homogeneous medium with wind).}
Suppose $R_H$, $R_B$, $U>0$, $S$, $\Theta_w$, and $\Theta_s$
are spatially constant, with the preferred direction defined and the
stated positivity and strong-convexity conditions on the full
angular rate $r(\theta)$ satisfied. Then the
coefficients $a,b,c,g,\lambda$ and the direction $\Theta$ are constant,
and Eq.~\eqref{model4wavelet} reduces to
\begin{equation}\label{fram2}
r(\theta)=\frac{c a^2\cos\delta+ab\sqrt{a^2\cos^2\delta+b^2\sin^2\delta-c^2\sin^2\delta}}{a^2\cos^2\delta+b^2\sin^2\delta}
+g\cos\delta,\qquad \delta=\theta-\Theta.
\end{equation}
The indicatrix is independent of $x$, and the translation-invariant Finsler metric is
\[
F(v)=\frac{\lVert v\rVert}{r(\arg v)}\quad(v\ne0),\qquad F(0)=0.
\]
The outward Finsler-unit propagation velocity is determined by
Eq.~\eqref{model3-support}, with $\mathcal I_x$ independent of $x$.
The corresponding fire ray is $\gamma_x(t)=x+t v_x$, for $t\ge0$.
The burned region at time $\tau$ is
\begin{equation}\label{model2front}
\begin{aligned}
K_\tau=K_0\cup\bigl\{x+t v_x:\;&x\in\partial K_0,\quad 0\le t\le\tau,\\
&\gamma_x|_{[0,t]}\text{ is minimizing from }K_0\bigr\}.
\end{aligned}
\end{equation}
The fire front at time $\tau$ is $\partial K_\tau$. The minimizing
condition excludes ray continuations after they cease to realize
the forward Finsler distance from the initial burned region.

\medskip\noindent\textit{Proof.}
Since the coefficients and preferred direction are independent of
position, so is $F$. Its geodesics are straight lines. For a smooth
convex initial burned region, all outward Finsler-normal rays remain
minimizing from $K_0$. For a nonconvex initial region, only their
minimizing segments contribute to Eq.~\eqref{model2front}, and the
fire front is the boundary of the reached region rather than the set
of unrestricted ray endpoints. The supporting-point condition and
its unit-speed and outward Finsler-normal consequences are precisely
those of the general construction, now with no base-point dependence.
\hfill$\square$

\medskip\noindent\textbf{Corollary 2 (heterogeneous medium without wind).}
Suppose $U(x)=0$ and $S(x)>0$ throughout the modeled region. Let
$R_H(x)>0$ and $R_B(x)>0$ be $C^2$ functions satisfying
$R_H(x)/3<R_B(x)\le R_H(x)$. Assume that the directional field
$g(x)(\cos\Theta_s(x),\sin\Theta_s(x))$ is $C^2$ (even where
$g(x)=0$), with $\Theta_s$ defined wherever this field is nonzero.
At points where $g=0$, interpret $g\cos(\theta-\Theta_s)$ through
its smooth vector-field extension rather than evaluating an undefined
angle. Then $a(x)=b(x)$, $c(x)=0$, and
\begin{equation}\label{model3wavelet}
 r(x,\theta)=b(x)+g(x)\cos(\theta-\Theta_s(x)),
\quad
\overline{\mathcal{I}_x}=\{\rho e_\theta:0\le\rho\le  r(x,\theta),\
0\le\theta<2\pi\},
\end{equation}
where
$b(x)=[R_H(x)+R_B(x)]/2$ and
$g(x)=[R_H(x)-R_B(x)]/2$.

The initial data for the corresponding Finsler metric $F_x(v)=\frac{\lVert v\rVert}{ r(x,\arg v)}$ satisfy
\begin{equation}\label{alp2}F_x(v_x)=1,\end{equation}
\begin{equation}\label{alp3}
g_{x,v_x}(v_x,u)=0\quad(u\in T_x\partial K_0).
\end{equation}
At a fixed $x$, the angular equation reduces to
\begin{equation}\label{model3-normal-explicit}
\lVert v_x\rVert= r(x,\theta),\qquad
u\cdot\left[b(x)e_\theta+g(x)e_{2\theta-\Theta_s(x)}\right]=0.
\end{equation}
The outward solution is selected by Eq.~\eqref{model3-support}.
Fire rays are generally curved geodesics, and the exact reached front
is $\partial K_\tau$ from Eq.~\eqref{model3front}. Spatial derivatives enter the geodesic
spray and do not enter this pointwise angular condition. 

\medskip\noindent\textit{Proof.}
Since $U(x)=0$ and $S(x)>0$, the parameter definitions give
$U^*(x)=0$, $a(x)=b(x)$, $\lambda(x)=1$, $c(x)=0$,
$g(x)=[R_H(x)-R_B(x)]/2$, and $\Theta(x)=\Theta_s(x)$.
Substitution into Eq.~\eqref{model4wavelet} yields
$r(x,\theta)=b(x)+g(x)\cos(\theta-\Theta_s(x))$.
Thus the local indicatrices, Finsler metrics, outward normal initial
data, geodesic equations, and reached regions agree with the stated
no-wind formulation. The pointwise condition
$R_H(x)/3<R_B(x)\le R_H(x)$ ensures the required strong convexity.
\hfill$\square$

\medskip\noindent\textbf{Corollary 3 (homogeneous medium without wind).}
\label{cor:homogeneous-no-wind}
Let the wind speed vanish, $U=0$, and assume $S>0$. Then
$\lambda=1$, $a=b$, $c=0$, $\Theta=\Theta_s$, and
$g=(R_H-R_B)/2$. The directional propagation rate and the
corresponding indicatrix at each point $x$ reduce to
\begin{equation}\label{ind}
r(\theta)=b+g\cos(\theta-\Theta_s),
\qquad
\mathcal{I}_x=\{r(\theta)e_\theta:0\le\theta<2\pi\}.
\end{equation}
If $3R_B>R_H$ and $R_H\ge R_B>0$, the resulting indicatrix is
strongly convex. Consequently, the Finsler metric, outward
Finsler-normal propagation velocities, and fire front evolution
coincide with those of the homogeneous no-wind model.

\begin{proof}
For $U=0$ and $S>0$, the parameter formulas yield
$U^*=0$, $a=b$, $\lambda=1$, $c=0$, $g=(R_H-R_B)/2$,
and $\Theta=\Theta_s$. Substitution into Eq.~\eqref{model4wavelet}
gives Eq.~\eqref{ind}. The translation-invariant Finsler metric
therefore takes the form
\begin{equation}\label{fins1}
F(v)=\frac{\lVert v\rVert}{r(\arg v)}
\quad(v\ne0),\qquad F(0)=0.
\end{equation}
The curvature numerator $r^2+2(r')^2-rr''$ attains its
minimum $(b-g)(b-2g)$ for $0\le g<b$. Under
$R_H\ge R_B>0$, the condition $g<b/2$ is equivalent to
$3R_B>R_H$ and guarantees strong convexity.
The outward Finsler-normal velocity satisfies
\begin{equation}\label{model1-normal-explicit}
\lVert v_x\rVert=r(\theta),\qquad
u\cdot[r(\theta)e_\theta-r'(\theta)e_\theta^\perp]=0,
\end{equation}
where $u$ is any nonzero tangent vector to $\partial K_0$
at $x$. Equivalently,
\begin{equation}\label{model1-normal-trig}
u\cdot\left[b e_\theta+g e_{2\theta-\Theta_s}\right]=0.
\end{equation}
The appropriate outward solution is selected by the supporting-point
maximization; an inverse-cosine expression alone generally does not
identify it. Since the indicatrix and metric agree with those of
the general homogeneous model in this limit, so do its normal
rays, their minimizing segments, and the resulting reached region.
\end{proof}

\begin{remark}\label{rem:exact}
The general model permits two equivalent descriptions of exact fire front propagation: tracing minimizing outward Finsler-normal geodesics from the initial front, or using the forward-distance function to construct point-source wavefronts and applying Huygens' principle. The latter construction accounts for geodesics in all Finsler-unit directions, whereas the former uses the initially outward normal directions. 

In a spatially homogeneous, time-independent medium, the Finsler metric is translation invariant and its geodesics are straight lines. Writing $\mathcal I=\{r(\theta)e_\theta:0\le\theta<2\pi\}$ for the unit indicatrix at point $x$ and $\overline{\mathcal I}$ for the filled body it bounds, the exact point-source reached region at time $T$ is $x+T\overline{\mathcal I}$, with wavefront $x+T\mathcal I$. Thus, the  indicatrix is the Huygens wavelet. This is the construction used in the numerical examples and applications under spatially uniform conditions; the polygonal representation approximates the wavelet geometry when discretized.

In heterogeneous media, the exact wavelets are instead determined by the spatially varying forward distance, and locally frozen indicatrices generally provide only a numerical approximation. If environmental conditions are updated between time intervals, the reached region from one interval serves as the initial region for the next. 
\end{remark}

\subsection{Geodesic Cut Times and Numerical Limitations}\label{sec}

For a time-independent, positive, strongly convex Finsler metric,
let $\gamma_s(t)$ be the outward unit-speed Finsler-normal geodesic
from $x_s\in\partial K_0$. Define
\[
E_t=\{\gamma_s(t):x_s\in\partial K_0\},\qquad
K_t=\{x\in M:d_F(K_0,x)\le t\},\label{reachedgeneral}
\]
where $d_F$ is the forward Finsler distance. While these geodesics
remain minimizing, their endpoints describe the advancing front.
Beyond a cut time, however, a ray may cease to minimize the distance
from $K_0$, for example when minimizing rays meet or at a focal point
of the normal family. A cut point may still lie on the front, whereas
later endpoints of the same ray need not. Thus, the fire front is
$\partial K_t$, which need not coincide with the unrestricted endpoint
set $E_t$. This distinction is important for nonconvex initial fronts
and heterogeneous media.

Our numerical scheme approximates reached regions by taking unions of
sampled wavelet polygons and extracting their boundaries; it does not
integrate geodesics. Under homogeneous, time-independent conditions,
the  indicatrix is the exact point-source wavefront and its filled
body is the exact wavelet, apart from polygonal discretization. In
heterogeneous media, locally frozen wavelets also approximate the
spatially varying propagation. If environmental inputs change between
intervals, each updated reached region becomes the initial region for
the next interval.

Further investigation of the scheme's applicability to complex topological configurations, including non-burnable barriers, persistent interior holes, and disconnected ignitions caused by spotting, represents a promising direction for future research. 

\subsection{Flowchart}\label{flowchart}

The following flowchart summarizes the main steps of the proposed wildfire front
propagation framework. Starting from prescribed spread rates and environmental
inputs, the directional wavelet is constructed and its Finsler admissibility
conditions are examined. The advancing front can then be characterized either
by minimizing outward Finsler-normal geodesics or by the boundary of the
forward-distance reached region through Huygens' principle.

In homogeneous, time-independent media, the local indicatrices provide exact
point-source wavelets. In heterogeneous media, the exact reached region is
determined by the position-dependent forward distance, while the numerical
implementation approximates it using unions of sampled locally frozen wavelets.
The resulting front is updated at each propagation step, and environmental
parameters may be revised between successive time-frozen intervals.

\begin{center}
\scalebox{0.83}{%
\begin{tikzpicture}[
  every node/.append style={text=blue,font=\small,align=center},
  bx/.style={process,text width=7.3cm,inner sep=5pt},
  dc/.style={decision,text width=3.5cm,inner sep=2pt},
  branch/.style={process,text width=5.1cm,inner sep=5pt}
]
\node (start) [startstop] at (0,0)
{1. Start: given initial burned region $K_0$
(front $\partial K_0$)};

\node (inputs) [bx] at (0,-1.8)
{2. Supply $R_H$, $k$ (with $R_B=kR_H$),
wind, slope, and directional data;
compute $a,b,c,g,\lambda,\Theta$
and the local indicatrix};

\node (admiss) [dc] at (0,-4.4)
{3. Is the directional wavelet positive
and strongly convex?};

\node (invalid) [startstop,text width=4.2cm]
at (5.3,-4.4)
{Stop: inadmissible inputs or wavelet.
Revise the parameters or formulation
before restarting};

\node (hom) [dc] at (0,-7.0)
{4. Spatially homogeneous conditions?};

\node (homcase) [branch] at (-4.3,-9.7)
{5. Position-independent Finsler metric;
geodesics are straight and point-source
wavelets are exact indicatrices};

\node (hetcase) [branch] at (4.3,-9.7)
{5. Position-dependent Finsler metric;
geodesics are generally curved;
frozen local wavelets approximate
forward-distance wavelets};

\node (method) [dc] at (0,-12.4)
{6. Trace fire rays?};

\node (rays) [branch] at (-4.3,-15.3)
{7. Find outward Finsler-unit normal
$v_x$ on $\partial K_0$ by the support
condition (Eq.\ref{model3-support}); solve Eq.~\eqref{geo-fin}
(straight rays in homogeneous media)};

\node (wave) [branch] at (4.3,-15.3)
{7. Apply Huygens' principle using
forward-distance reached regions;
in numerical propagation, use locally frozen
heterogeneous wavelets};

\node (rayfront) [branch] at (-4.3,-18.5)
{8. Retain minimizing ray segments
up to cut times; determine the front
as the boundary of the reached region};

\node (wavefront) [branch] at (4.3,-18.5)
{8. Extract the boundary of the reached
region; polygon unions approximate
the numerical update};

\node (cont) [dc] at (0,-21.2)
{9. Continue propagation?};

\node (update) [bx] at (0,-23.8)
{10. Replace $K_0$ by the updated
burned region. If environmental inputs
change, recompute the coefficients and
metric for the next time-frozen interval};

\node (stox) [startstop] at (4.5,-21.2)
{Stop};

\draw [arrow] (start) -- (inputs);
\draw [arrow] (inputs) -- (admiss);

\draw [arrow] (admiss) --
node[right] {yes} (hom);

\draw [arrow] (admiss.east) --
node[above] {no} (invalid.west);

\draw [arrow] (hom) -|
node[pos=0.2,above] {yes} (homcase);

\draw [arrow] (hom) -|
node[pos=0.2,above] {no} (hetcase);

\draw [arrow] (homcase.south) --
++(0,-0.45) -| (method.west);

\draw [arrow] (hetcase.south) --
++(0,-0.45) -| (method.east);

\draw [arrow] (method) -|
node[pos=0.2,above] {yes} (rays);

\draw [arrow] (method) -|
node[pos=0.2,above] {no} (wave);

\draw [arrow] (rays) -- (rayfront);
\draw [arrow] (wave) -- (wavefront);

\draw [arrow] (rayfront.south) --
++(0,-0.5) -| (cont.west);

\draw [arrow] (wavefront.south) --
++(0,-0.5) -| (cont.east);

\draw [arrow] (cont) --
node[right] {yes} (update);

\draw [arrow] (cont) --
node[above] {no} (stox);

\draw [arrow] (update.west) --
++(-6.6,0) |- (inputs.west);

\end{tikzpicture}%
}
\end{center}


\section{Illustrative examples}\label{examp}
The three examples illustrate the Finsler construction under prescribed
synthetic inputs; they are not calibrated wildfire forecasts. Lengths are
in km, times in hours, speeds in km/h, and angles are measured
counterclockwise from East.

\subsection{Prescribed rate inputs for the illustrative scenarios}\label{rateconversion}
The head and backing rates $R_H$ and $R_B$ specify the local directional
wavelets. The inputs are synthetic and are not estimated from observations.
\begin{table}[htbp]
\centering\small
\begin{tabular}{@{}p{0.25\linewidth}p{0.29\linewidth}p{0.29\linewidth}r@{}}
\toprule
Scenario & $R_H$ (km/h) & $R_B$ (km/h) & $U$ (km/h)\\\midrule
Ex.~\ref{ex1}(a) & $1.45[1.8-0.6\cos(x+y)]$ & $0.55[1.8-0.6\cos(x+y)]$ & $0$\\
Ex.~\ref{ex1}(b) & $1.30[3.5+\cos^2 y]$ & $0.70[3.5+\cos^2 y]$ & $0$\\
Ex.~\ref{ex1}(c), first 10 h & $1.30[1.8-\cos(x+y)]$ & $0.70[1.8-\cos(x+y)]$ & $0$\\
Ex.~\ref{ex1}(c), next 20 h & $5.280$ & $0.64627$ & $7$\\\bottomrule
\end{tabular}
\caption{Head and backing spread rates for the illustrative scenarios. }\label{rateinputs}
\end{table}

\begin{example}
    \label{ex1}
The initial burning region is a circle of radius $0.5$ km.
Figure~\ref{slope} illustrates three scenarios. 
In panel (c), the first $10$ h is without wind and $\Theta_s=\pi/4$
(northeast). 
For the following $20$ h, 
$\Theta_w=3\pi/4$ (northwest) and $\Theta_s=\pi/4$
(northeast). The rates follow from the illustrative input-rate calculation in Section~\ref{rateconversion}.

\begin{figure}[h!]
\centering
\includegraphics[width=\textwidth]{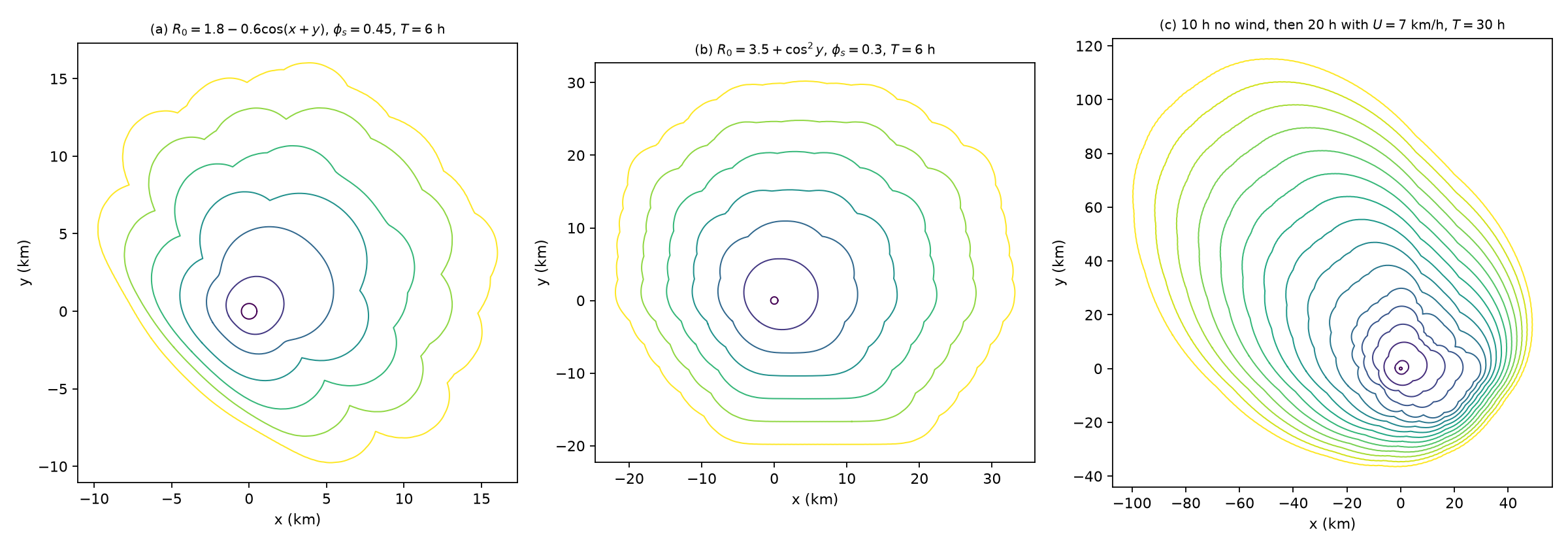}
\label{slope}
\end{figure}

\end{example}

**Example 4.4. Wildfire propagation over a real heterogeneous landscape in northwestern Ontario, Canada.**

To investigate the applicability of the proposed geometric framework to real topographical conditions, we consider a geographical region in northwestern Ontario, Canada, centred at latitude 51.0000° N and longitude 91.0000° W. The numerical experiment is motivated by wildfire activity in northwestern Ontario during July 2012, with particular attention to the period from 13 to 15 July.

The computational domain covers an area of 20 km × 20 km and is discretized into 1,600 square cells of side length 500 m. Terrain elevation data are obtained from the Medium Resolution Digital Elevation Model (MRDEM-30) provided by Natural Resources Canada, with a spatial resolution of approximately 30 m.

The elevation within the selected region varies between 362.98 m and 432.89 m, with a mean elevation of 378.33 m. To determine the local terrain inclination, a plane is fitted by least squares to the elevation samples within each computational cell. The corresponding terrain gradient is used to determine the local inclination factor and uphill direction according to the formulation introduced in Section 2.

The resulting inclination angles range from approximately 0° to 5.24°, with a mean value of 0.60°. Although the selected landscape is predominantly gently inclined, the spatial variability of its terrain gradients provides heterogeneous directional conditions for wildfire propagation.

Unlike the preceding illustrative examples, the terrain inclination magnitude and uphill direction are obtained from actual geographical elevation data rather than prescribed synthetic functions. The environmental coefficients associated with the terrain remain fixed throughout the simulation, whereas the wind conditions may be updated between successive propagation intervals.

The purpose of this numerical experiment is to illustrate the applicability of the proposed directional wavelet construction to a real heterogeneous landscape. It should be distinguished from an empirical reconstruction or validation of a historical wildfire event, which would require independently established initial conditions, spread rates, fuel characteristics, and observed fire perimeters.

\subsection{\textcolor{black}{Numerical algorithm}}\label{algsec}
\rev{The reached region $K_t$ and the fire perimeter $\partial K_t$ are advanced as follows.
\begin{enumerate}
\item \textbf{Initial perimeter.} The initial burning set $K_0$ is a polygon (a circle approximated by 360 vertices, a lobed curve $\rho=1+0.45\cos3\varphi$, or a grid cell). Its boundary is resampled to $N=240$ equally spaced vertices $x_1,\dots,x_N$.
\item \textbf{Local wavelet.} At each sampled front point and time step,
supply $R_H$, $R_B$, wind speed $U$, wind direction $\Theta_w$, terrain
slope angle $\phi$, and uphill direction $\Theta_s$.
Use the formulation in Appendix~\ref{apen3} to compute the directional
wavelet. For wind cases, set $U^*=U/(3.6\,\mathrm{km/h})$
and use the phenomenological wind--slope weight $\mu=1$.
No slope or wind enhancement factors are required by the present model;
$R_H$ and $R_B$ must be supplied separately. Check their positivity,
the preferred-direction definition, and angular convexity conditions
before each propagation update. Finite angular sampling provides only
an admissibility screen, not a proof over every angle.
Sample the directional wavelet at $160$ equally spaced angles to obtain
the local polygon $x_i+\Delta t\,\overline{\mathcal I_{x_i}}$.
\item \textbf{Time step.} $\Delta t=0.5$ h in Example~\ref{ex1}(a,b) and Example~\ref{ex2}, $\Delta t=1$ h in Example~\ref{ex1}(c) and Example~\ref{proof} (fronts are plotted at the intervals stated in the captions). \rev{The step is chosen so that the displacement $\Delta t\,R_H$ stays small compared with the length scale of the front and of the coefficient variation; the effect of $\Delta t$ is examined empirically for one scenario in Section~\ref{convsec}. For a homogeneous medium the scheme was also compared with the exact  wavelet (symmetric-difference area $0.3\%$ of the area).}
\item \textbf{Perimeter update.} $K_{t+\Delta t}=K_t\cup\bigcup_i (x_i+\Delta t\,\overline{\mathcal I_{x_i}})$ (polygon union). The new perimeter is the exterior boundary of $K_{t+\Delta t}$; interior holes (unburnt islands) are not propagated. \rev{The reached region is retained between steps (each update is the union of the region and the new wavelets), but the next set of front points is taken from the exterior boundary of the largest connected component only. Interior holes may therefore be filled by overlapping wavelets and are never used as propagation sources; if the union becomes a multipolygon, all components are kept in the region but only the largest propagates. The implemented update is therefore a \emph{single-exterior-front approximation}: it does not propagate every connected component, and holes are neither propagated nor preserved as separate fronts. Disconnected burning regions are thus not treated by Huygens propagation.} Steps 2--4 are repeated, switching the parameters at the phase times.
\item \textbf{Grid example.} For a square cell of side $0.5$ km, a $4\times4$ array of equally spaced sample points is placed in the cell, and the cell is marked as burnt if at least $8$ of the $16$ points lie in $K_t$. This is an approximate $50\%$ area threshold based on spatial sampling, not an exact area-fraction calculation. For comparison, the exact polygon--cell intersection area was also computed at $t=21$ h: it gives $484$ burnt cells (area fraction $\ge50\%$) against $485$ from the sampling test.
\end{enumerate}}

\subsection{\textcolor{black}{Numerical refinement study}}\label{convsec}
\rev{The numerical consistency of the front-evolution algorithm under refinement is assessed for the no-wind scenario of Example~\ref{ex1}(a): $R_H=1.45[1.8-0.6\cos(x+y)]$ km/h, $R_B=0.55[1.8-0.6\cos(x+y)]$ km/h, an initial circle of radius $0.5$ km and a final time $T=6$ h. The time step was $\Delta t=2,1,0.5,0.25$ h with the spatial discretization fixed ($N=240$ front vertices, $160$ wavelet angles). The reference solution used $\Delta t=0.125$ h with $N=480$ and $320$ angles. For each run, Table~\ref{conv} reports the final burnt area, the perimeter length of the final polygon, and the symmetric-difference area between the final region and the reference region, divided by the reference area. A last row repeats $\Delta t=0.25$ h with the finer spatial discretization.}

\begin{table}[h!]
\centering
\rev{\small\begin{tabular}{@{}lcccc@{}}
\hline
$\Delta t$ (h) & $N$, angles & area (km$^2$) & perimeter (km) & sym. diff./ref.\ area \\
\hline
2 & 240, 160 & 435.5 & 78.8 & 15.3\% \\
1 & 240, 160 & 465.2 & 81.6 & 8.7\% \\
0.5 & 240, 160 & 441.9 & 81.2 & 2.0\% \\
0.25 & 240, 160 & 436.7 & 83.5 & 0.87\% \\
0.25 & 480, 320 & 442.1 & 82.0 & 1.1\% \\
\hline
reference: 0.125 & 480, 320 & 438.1 & 84.1 & --- \\
\hline
\end{tabular}}
\caption{\rev{Time-step refinement and numerical consistency for the heterogeneous scenario of Example~\ref{ex1}(a), $T=6$ h; differences are relative to a numerical reference solution, not to an exact solution. The symmetric difference (``sym. diff.'') is the area of the set where the final region and the reference region differ. It decreases as $\Delta t$ is refined from $1$ to $0.25$ h; the burnt area alone is not monotone. The perimeter length of a polygon depends on the vertex resolution and is shown for information only.}}\label{conv}
\end{table}

\rev{With $\Delta t=1$ h the final region differs from the reference by about $9\%$ of its area, with $\Delta t=0.5$ h by $2\%$ and with $\Delta t=0.25$ h by $0.87\%$. These values describe the empirical refinement behavior relative to a numerical reference solution, not an error bound. The remaining difference of about $1\%$ reflects the interacting effects of the time step, perimeter sampling, angular sampling and the discretization of the reference itself, and is not attributed to a single source; indeed, at $\Delta t=0.25$ h the finer spatial resolution gave a slightly larger symmetric difference ($1.1\%$ versus $0.87\%$), so no further monotone improvement is observed at this level. The study covers one scenario only, no convergence theorem is claimed, and the time-step study was not repeated for the other examples. As a homogeneous-medium check only, the scheme reproduces the exact  wavelet with a symmetric-difference area of $0.3\%$. The qualitative conclusions (shape, orientation and elongation of the fronts) are unchanged by the refinement.}

\subsection{\textcolor{black}{Computational resources and reproducibility}}\label{ressec}
The numerical construction uses sampled wavelets and polygon-union operations,
which can be implemented with NumPy, Matplotlib, and Shapely. The inputs and
discretization settings are specified above. Reproducibility of the numerical
figures and tabulated results requires execution of the implementation with
these inputs and verification of the admissibility conditions.

\rev{Examples \ref{ex1} and \ref{ex2} illustrate the geometry of the proposed wavefronts under prescribed homogeneous and spatially varying coefficients, with and without wind. They do not estimate or compare fire intensity.}

In Example \ref{proof}, we discretize the terrain into a grid to track the progression of the fire. Using the flowcharts presented in \ref{flowchart}, we model the fire’s spread and apply this model to the grid to simulate its advance.
\subsection{Simulation results and geometric illustrations}

\begin{example}\label{ex2}
The initial front is the nonconvex three-lobed curve
$\rho=1+0.45\cos3\varphi$ km. The homogeneous scenario uses
$R_H=6.020$ km/h, $R_B=0.737$ km/h, and $U=7$ km/h,
with wind and uphill directions aligned at
$\Theta_w=\Theta_s=\pi/4$. Figure~\ref{slope2} compares the Huygens reached-region boundary with outward Finsler-normal straight rays. Figure~\ref{slope2old} additionally presents the two constructions in separate panels at $\Delta t=2$ h.
\begin{figure}[h!]
\centering
\includegraphics[width=\textwidth]{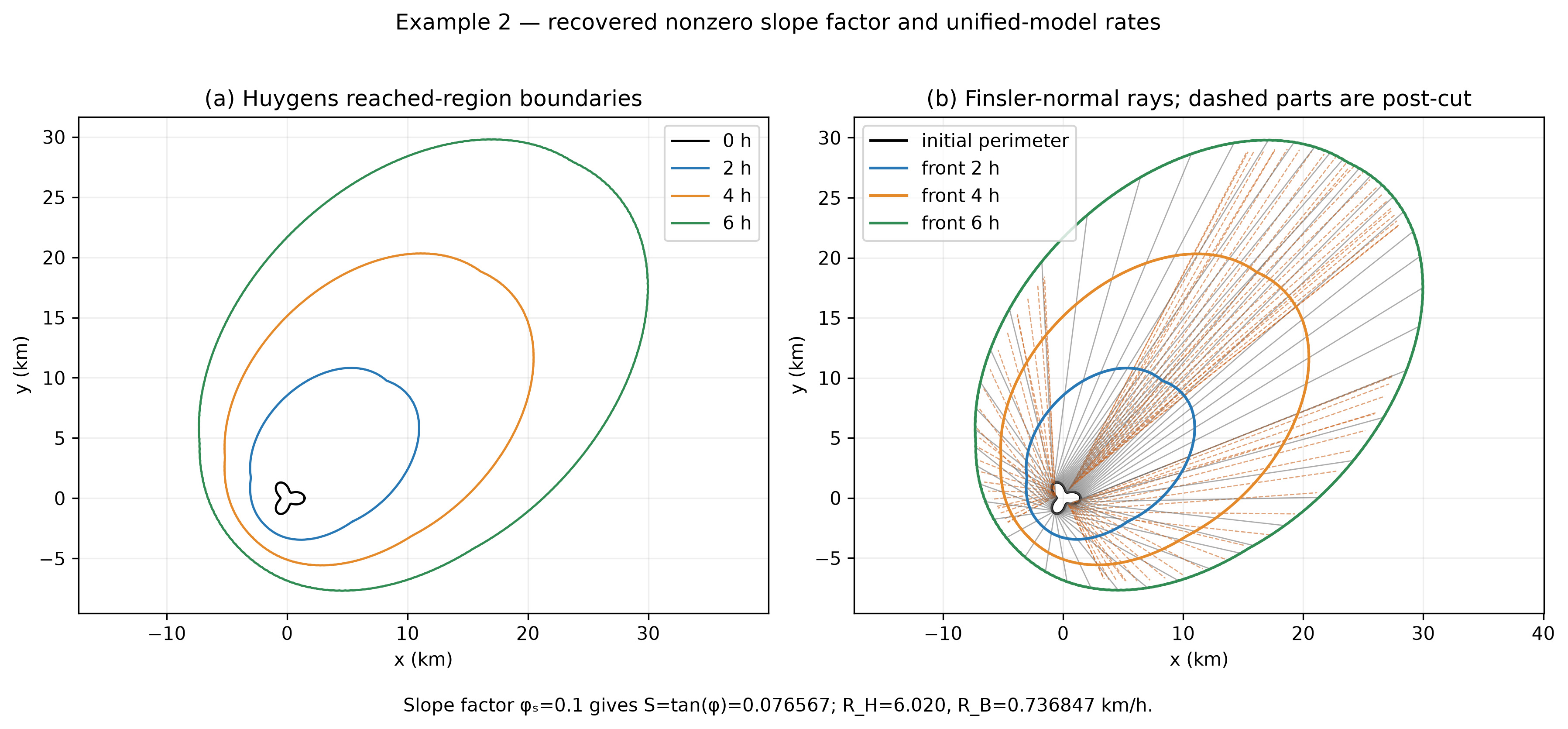}
\caption{{Example~\ref{ex2}: geometric comparison of the
Huygens reached-region boundary and straight fire rays at $t=2,4,6$ h.
Prescribed $R_H=6.020$ km/h, $R_B=0.737$ km/h, $U=7$ km/h,
$\Theta_w=\Theta_s=\pi/4$; initial front
$\rho=1+0.45\cos3\varphi$ km. The rates are prescribed synthetic inputs.}}\label{slope2}
\end{figure}

\begin{figure}[h!]
			\centering
           \begin{subfigure}[b]{0.45\textwidth}
         \centering
         \includegraphics[width=.65\textwidth]{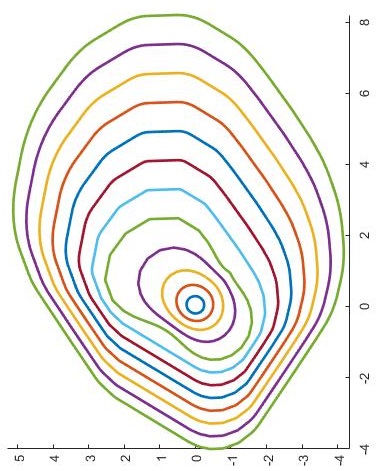}
         \subcaption{Fire model applying Huygens' principle.}\label{r1}
     \end{subfigure}
     \hfill
     \begin{subfigure}[b]{0.45\textwidth}
         \centering
         \includegraphics[width=.65\textwidth]{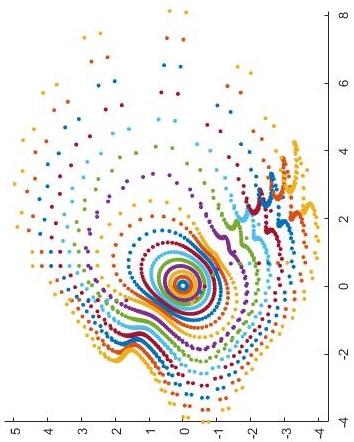}
         \subcaption{Fire model by finding the fire rays.}\label{r2}
     \end{subfigure}
           \caption{Comparing the methods of modeling using Huygens' principle and rays with $\Delta t=2$ hours (Example \ref{ex2}). }\label{slope2old}
     \end{figure}

\end{example}

The grid example applies the same polygon-union update described in Section~\ref{algsec}.

\begin{example}[\textbf{Proof of Concept}]\label{proof}

We simulated a simplified scenario of fire spread over a discretized landscape. The terrain is divided into grid cells, and the fire's propagation is modeled using data on topography, vegetation characteristics for the western-grass fuel type following \cite{andrews2018rothermel},
and wind conditions that are prescribed illustrative scenario parameters (see below). Fires are assumed to originate from a single grid cell within the landscape.

\rev{In Fig.~\ref{cells}, a cell is highlighted in red if at least $8$ of its $16$ equally spaced sample points ($4\times4$ array) lie in the reached region; this approximates, but is not an exact test of, a $50\%$ affected-area threshold.} \textcolor{black}{Cells that ignite during the current interval are shown in red, while those ignited in earlier intervals are shown in light red.} The sequence of plots in Fig.~\ref{cells} shows the simulated fire spread \rev{at the end of each of the four phases ($t=5,12,16,21$ h), illustrating how grid-based modeling can track fire progression over time}. 

The simulation parameters are as follows: the vegetation type is classified as western grass, and wind speed and direction vary across four distinct phases over a \rev{21-hour} period \rev{(5, 7, 4 and 5 h)}. Here, $\Theta_s$ represents the uphill direction, and $\Theta_w$ denotes the wind direction \rev{measured counterclockwise from the $x$-axis (East)}. The specific parameters used for each propagation phase are summarized in Table~\ref{param}. \rev{The fuel parameters follow the western-grass values of \cite{andrews2018rothermel}; the wind speeds and directions are prescribed illustrative scenario parameters selected to demonstrate the behavior of the proposed geometric construction. They are not observationally derived measurements and are not taken from a site-specific dataset.}

\rev{All four phases use the prescribed synthetic head and backing
rates, wind speeds, and directional data in Table~\ref{param}.
The ratios $R_B/R_H$ are approximately $0.0790$, $0.0646$,
$0.0503$ and $0.0450$. Finite-grid angular checks provide an admissibility screen,
not continuous-domain convexity certification. The fire begins in one cell of side $0.5$ km
and the finite-step perimeter update uses $\Delta t=1$ h.
The cell counts used for the numerical comparison are
$1,15,136,290,485$ at $t=0,5,12,16,21$ h.
No satellite calibration is claimed.}

\begin{table}[h!]
\centering
\begin{tabular}{lcccc}
\hline
\textbf{Parameter} & \textbf{phase 1} & \textbf{phase 2} & \textbf{phase 3} & \textbf{phase 4} \\
\hline
$R_B (km/h)$       & \textcolor{black}{0.032} & \textcolor{black}{0.052} & \textcolor{black}{0.052} & \textcolor{black}{0.038}      \\
$\Theta_s (rad)$    & 0.2618 & 0.1745 & 0.3491 & 0.1745 \\
$R_H (km/h)$       & \textcolor{black}{0.405} & \textcolor{black}{0.805} & \textcolor{black}{1.034} & \textcolor{black}{0.844} \\
$U (km/h)$         & 12     & 15     & 17     & 20     \\
$\Theta_w (rad)$  & 1.50   & 1.20   & 1.20   & 1.30   \\
$T (h)$ time         & 5      & 7      & 4      & 5      \\
\hline 
\end{tabular}
\caption{\textcolor{black}{Parameters used for four phases of propagation (total 21 h); $R_H$ and $R_B$ prescribed for this synthetic example (with $k=R_B/R_H$); angles counterclockwise from East.}}
\label{param}
\end{table}
\begin{figure}[h!]
\centering
\includegraphics[width=\textwidth]{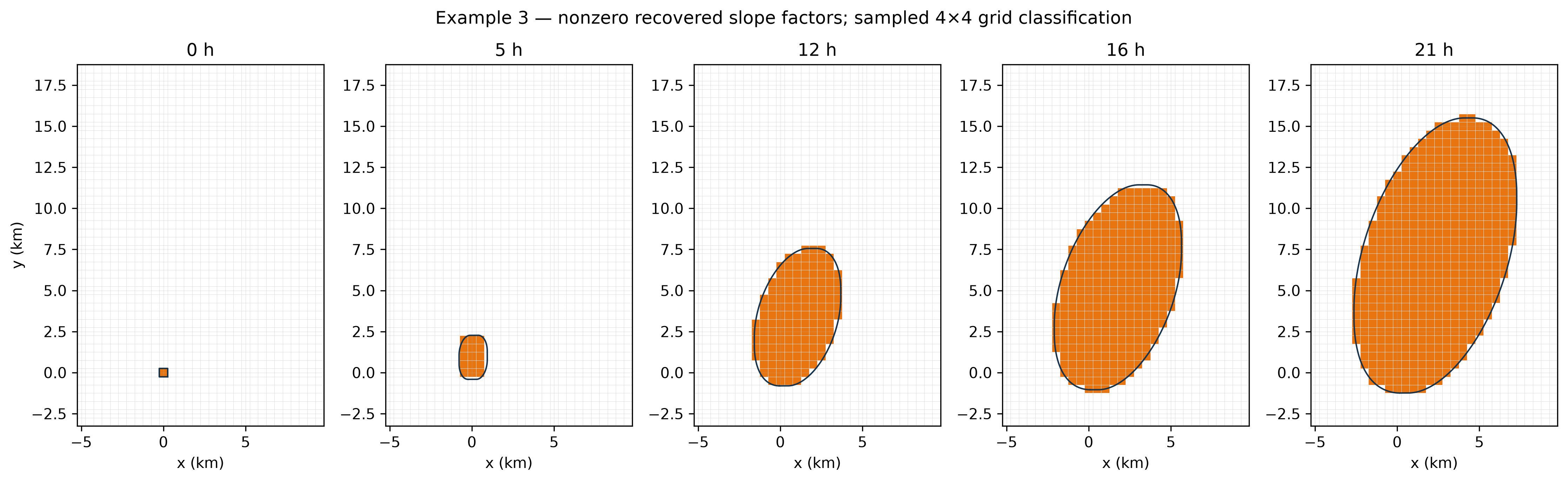}
\caption{\textcolor{black}{Example~\ref{proof}: grid propagation with cells of $0.5$ km (axes in km), parameters of Table~\ref{param}. Panel (a): ignition cell; (b)--(e): state at the end of phases 1--4 ($t=5,12,16,21$ h); red cells ignited during the last phase, light-red cells earlier. Panel (f): fire perimeters every 1 h from $0$ to $21$ h (colour scale from dark to light). A cell is classified as burnt if at least 8 of its 16 equally spaced sample points ($4\times4$) lie inside the reached region; this approximates a 50\% burned-area threshold.}}\label{cells}
\end{figure}

\end{example}
\subsection{Engineering interpretation and potential applications}

\rev{From an engineering perspective, the proposed framework provides a
mathematical tool for examining how prescribed environmental conditions
affect the geometry of an advancing wildfire front. The directional
wavelet describes the local propagation response, while the Huygens-type
construction combines these local responses to obtain an evolving burned
region. This separation is useful because changes in the propagation law
can be studied without modifying the general front-evolution procedure.}

\rev{The four model configurations allow the effects of wind, slope, and
spatial heterogeneity to be examined individually or in combination.
Under homogeneous conditions, the resulting fronts provide reference
configurations for understanding the influence of directional anisotropy.
Under heterogeneous conditions, the local propagation parameters vary
across the domain, allowing the front to respond to prescribed changes
in the environment. Such experiments can help identify which assumptions
and parameter variations have the greatest influence on the computed
perimeter.}

\rev{The formulation may also provide a geometric component for future
simulation and decision-support studies. For example, computed burned
regions could be used as inputs to mathematical models concerned with
exposure assessment, spatial planning, or the evaluation of hypothetical
protection strategies. These applications would require additional
modeling, calibration, and validation and are not demonstrated by the
present numerical experiments.}

\rev{The current framework should consequently be regarded as an
idealized engineering-mathematics model rather than an operational
wildfire prediction system. Its principal value is the explicit
relationship between directional spread-rate parameters, admissible
propagation geometry, and the resulting front evolution. Validation
against observed fire perimeters and comparisons with established
simulators are necessary before conclusions about predictive performance
or practical decision-making can be drawn.}

\section{Conclusion and Final Remarks}\label{concl}
\rev{This study develops an explicit Finsler-geometric framework for wildfire front propagation on locally inclined landscapes. It connects prescribed head-fire and backing-fire rates to a two-dimensional directional wavelet across homogeneous and heterogeneous environments, with and without wind. The wavelet combines an offset-ellipse radial profile and a slope-driven cosine term through a vector wind--slope weighting rule. By construction, it recovers the specified rates $R_H$ and $R_B=kR_H$ along the preferred forward direction and its opposite, and reduces to the windless slope-driven and flat-terrain limits.}

\rev{The geometric formulation requires strict positivity and strong convexity of the local indicatrix. We establish the exact windless cosine-wavelet threshold $R_H<3R_B$ (equivalently, $k>1/3$) and formulate the full angular curvature condition for combined wind and slope, with numerical angular screening for the illustrative windy cases. The supporting-point maximization over the filled indicatrix uniquely determines the outward Finsler-normal velocity without directly solving a nonlinear orthogonality equation. In homogeneous media, geodesics are straight; in heterogeneous media, rays are generally curved, and the exact front is described by minimizing geodesics up to their cut times or by the boundary of the forward-distance reached region. Numerical illustrations with nonconvex initial perimeters, a refinement study, and multi-phase grid propagation demonstrate the geometric behavior of the resulting fronts.}

\rev{The framework is a mathematical rate-to-wavelet closure, not an empirical law for the spread rates. The rates $R_H$ and $R_B=kR_H$, together with wind and terrain data, are externally prescribed, and changes in environmental conditions are represented by piecewise-static updates. In the present phenomenological construction, changing the relative wind--slope orientation rotates the preferred propagation direction without directly changing the prescribed head and backing magnitudes. The numerical examples use synthetic inputs and are neither calibrated to observed wildfire perimeters nor benchmarked against operational simulators; consequently, they illustrate geometric consistency rather than operational forecasting accuracy.}

\rev{Future work can incorporate angle-dependent spread-rate laws that account for wind--slope interactions, and investigate cone-based Lorentz--Finsler formulations when backing spread vanishes. Further computational studies should examine refinement on heterogeneous landscapes, non-burnable barriers and front merging, and the integration of spread-rate estimates from remote-sensing observations. Field validation and systematic comparison with established wildfire simulators are also needed to assess predictive performance.}


	\section*{Acknowledgments}
	
	The  author gratefully acknowledges support from the São Paulo Research Foundation (FAPESP, Brazil), Process No. 2022/15371-3, Process No. 2024/00923-6, and 2026/01122-2.




\appendix
\renewcommand{\theequation}{\thesection.\arabic{equation}}
\section{Rate inputs and the directional formulation}\label{apen3}

The general model takes the local head-fire spread rate $R_H(x,t)>0$ as an externally supplied input, potentially estimated from satellite-derived fire front displacement or field observations. Wind speed $U(x,t)$ and direction $\Theta_w(x,t)$ are likewise supplied. The backing-fire rate is prescribed by
\[
R_B(x,t)=k(x,t)R_H(x,t),\qquad 0<k(x,t)\le1.
\]
The ratio $k$ is selected using vegetation information, empirical evidence, or stated expert assumptions; it is neither estimated nor validated in this study. Thus, the mathematical model starts from prescribed rates and environmental inputs rather than deriving $R_H$ or $R_B$ from a combustion or Rothermel-type rate law. A satellite-derived $R_H$ represents an estimated displacement rate over a specified spatial support and observation interval; uncertainties in that estimate are outside the present analysis.

The terrain slope entering the parameter $\lambda$ is the dimensionless gradient magnitude. For an elevation field $z=h(x,y)$,
\[
S=\|\nabla h\|=\sqrt{h_x^2+h_y^2}=\tan\phi,
\]
where $\phi$ is the slope angle relative to the horizontal, and $\Theta_s$ denotes the uphill direction wherever the gradient is nonzero. The definitions of $U^*$, $a,b,c,g,\lambda$, and the resultant direction $\Theta$ are given in the general model; they are not repeated here.

For implementation, angles follow the mathematical convention: counterclockwise from East, pointing in the direction of propagation. A meteorological wind direction $\Theta_{\rm met}$, measured clockwise from North and indicating where the wind comes \emph{from}, is converted to the downwind propagation direction by
\[
\Theta_w=(270^\circ-\Theta_{\rm met})\bmod 360^\circ.
\]
The wind--slope resultant is
\begin{equation}\label{prefdir}
(q_x,q_y)=(1-\lambda)(\cos\Theta_w,\sin\Theta_w)
+\lambda(\cos\Theta_s,\sin\Theta_s),\qquad
\Theta=\operatorname{atan2}(q_y,q_x).
\end{equation}
Only directions with nonzero weights need to be evaluated. If the resultant vanishes, the formulation does not assign a preferred direction; no fallback is assumed. For zero wind and positive slope, $\lambda=1$ and $\Theta=\Theta_s$; for zero slope and positive wind, $\lambda=0$ and $\Theta=\Theta_w$. When both wind and slope vanish, neither direction is evaluated: the adopted isotropic convention imposes $k=1$ and $r(x,\theta)=R_H(x)=R_B(x)$, independently of $\theta$. This is an additional modeling convention, not a consequence of the undefined ratio defining $\lambda$.

The wind--slope weighting parameter $\mu$ and the numerical factor $0.25$ in the definition of $a$ are phenomenological shape choices, not estimated physical constants. The rates $R_H$ and $R_B$ are specified independently of the relative wind--slope angle. In particular, rotating the wind direction relative to the slope rotates the resulting wavelet through $\Theta$ but does not, by itself, alter its prescribed head and backing rates or its shape. This is a limitation of the formulation: actual wind--slope interactions may also modify spread rates. The case $k=0$ gives a vanishing backing rate and lies outside the positive strongly convex Finsler setting.

\begin{figure}[h!]
\centering
\includegraphics[width=\textwidth]{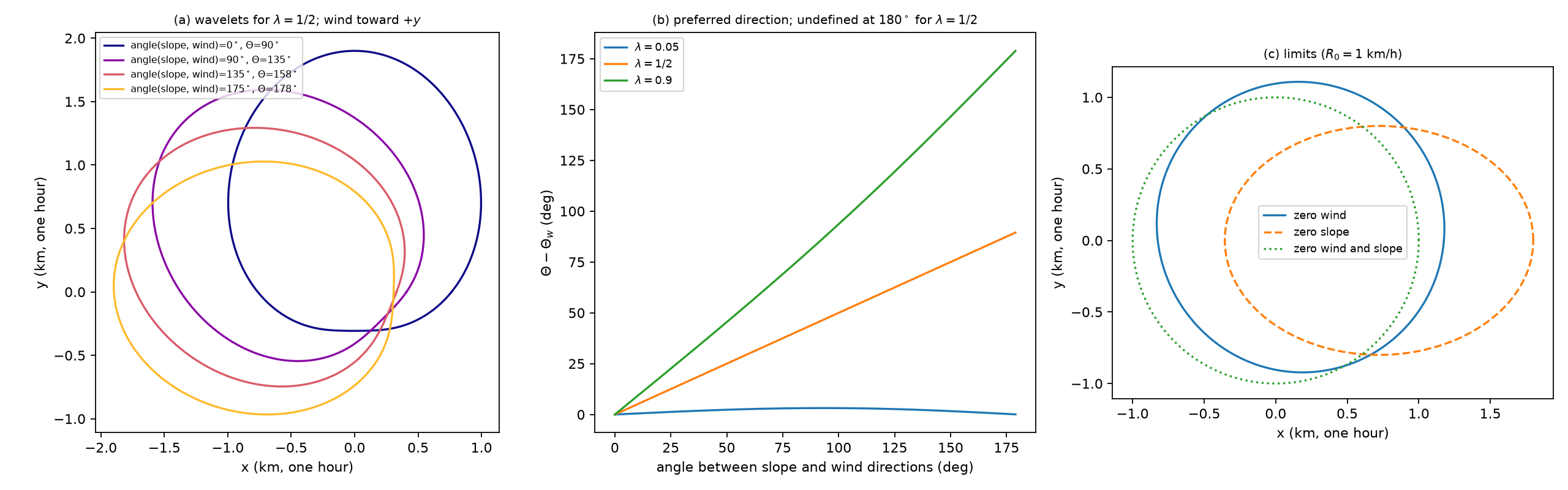}
\caption{Illustration of the wind--slope directional formulation. Panel (a) shows rigid rotations of prescribed synthetic wavelets with $U=5$ km/h and several wind--uphill angular separations; (b) shows the preferred-direction shift for selected weighting ratios $\lambda=0.05,0.5,0.9$; (c) illustrates no-wind, no-slope, and isotropic limits. When the wind--slope resultant vanishes, no preferred direction is defined. The plotted historical wavelets have not been independently regenerated from a documented rate-only data set.}\label{align}
\end{figure}

\paragraph{Limiting configurations.}
The zero-wind, positive-slope limit gives $a=b$, $c=0$, and $r=b+g\cos(\theta-\Theta_s)$; the zero-slope, positive-wind limit gives $g=0$ and an offset ellipse. The zero-wind, zero-slope case follows the isotropic convention above. For exactly opposing wind and uphill directions, the resultant in Eq.~\eqref{prefdir} vanishes when $\lambda=1/2$; on either side of that value its direction follows the more heavily weighted contribution. Near cancellation, small changes in the inputs can cause large changes in $\Theta$. Figure~\ref{align} illustrates these configurations. The illustrated limiting cases are algebraic properties of the formulation, not empirical validation of wildfire forecasts.

\section{Finsler fundamental tensor and admissibility conditions}\label{geodesics}

This appendix gives the fiber Hessian and the positivity and strong-convexity tests for the directional-speed metric of the general model. At a fixed position $x$ and within a time-frozen interval, write $r(\theta)=r(x,\theta)$ and $\psi(\theta)=1/r(\theta)$. Then
\[
F_x(v)=\|v\|\psi(\arg v),\qquad v\ne0.
\]
The fundamental tensor $g_{x,v}=\tfrac12D_v^2F_x^2$, expressed in the orthonormal radial--angular frame $(e_\theta,e_\theta^\perp)$, is
\[
[g_{x,v}]=
\begin{pmatrix}
\psi^2 & \psi\psi'\\
\psi\psi' & \psi^2+(\psi')^2+\psi\psi''
\end{pmatrix},
\]
where primes denote angular derivatives. It is positive definite if and only if $\psi>0$ and $\psi+\psi''>0$. Since
\[
\psi+\psi''=\frac{r^2+2(r')^2-rr''}{r^3},
\]
the admissibility conditions are
\begin{enumerate}
\item[(A1)] $r(\theta)>0$ for every $\theta$;
\item[(A2)] $K(\theta):=r(\theta)^2+2[r'(\theta)]^2-r(\theta)r''(\theta)>0$ for every $\theta$.\label{cuvature}
\end{enumerate}
The curvature of the indicatrix is $\kappa=K/(r^2+(r')^2)^{3/2}$. Positive backing speed alone does not guarantee (A2) for the combined radial law. The condition $c<b$ ensures positivity of the offset-ellipse component, but not necessarily positivity and strong convexity after the cosine term is added.

\paragraph{Special cases and numerical screening.}
For the no-wind limit $r=b+g\cos\delta$, with $0\le g<b$, the curvature numerator has minimum $(b-g)(b-2g)$. Thus the indicatrix is strongly convex precisely when $g<b/2$, equivalent to $R_B>R_H/3$ or $k>1/3$. For the pure-wind case $g=0$, the indicatrix is an ellipse containing the origin when $c<b$ and is strongly convex. For mixed wind and slope, no general analytical admissibility criterion is established here; the numerical checks reported below sample the angular domain and are not proofs of positivity at every angle.

The historical mixed-case values in Table~\ref{adm} were obtained using $4096$ equally spaced angles and spectral (FFT) derivatives. The quantity $\kappa b$ is dimensionless. These values have not been independently recomputed from the original simulation script, and no angular-resolution sensitivity or certified interval bound is reported. The admissibility conditions are unchanged by multiplying $r$ by a positive constant, so normalized checks depend on the ratios $a/b$, $c/b$, and $g/b$, rather than the absolute scale of $R_H$.

\begin{table}[htbp]
\centering\small
\caption{Admissibility information expressed through rate inputs and wind speed. Curvature minima for windy cases are historical numerical screening values; ``screened'' does not imply certified admissibility.}\label{adm}
\begin{tabular}{@{}p{4.1cm}ccc p{2cm}@{}}
\toprule
Scenario & $R_B/R_H$ & $U$ (km/h) & min $\kappa b$ & Assessment \\
\midrule
Ex.~\ref{ex1}(a) & $0.3793$ & $0$ & $0.331$ & analytic: pass \\
Ex.~\ref{ex1}(b), (c) stage 1 & $0.5385$ & $0$ & $0.816$ & analytic: pass \\
Ex.~\ref{ex1}(c), stage 2 & $0.1224$ & $7$ & $0.670$ & screened \\
Ex.~\ref{ex2} & $0.1224$ & $7$ & $0.670$ & screened \\
Grid phase 1 & $0.0790$ & $12$ & $0.543$ & screened \\
Grid phase 2 & $0.0646$ & $15$ & $0.488$ & screened \\
Grid phase 3 & $0.0503$ & $17$ & $0.456$ & screened \\
Grid phase 4 & $0.0450$ & $20$ & $0.417$ & screened \\
No wind illustrative limit & $2/3$ & $0$ & $0.938$ & analytic: pass \\
Isotropic limit & $1$ & $0$ & $1$ & analytic: pass \\
No-wind inadmissible example & $1/4$ & $0$ & $-1.25$ & analytic: fail \\
Degenerate backing rate & $0$ & --- & --- & excluded \\
\bottomrule
\end{tabular}
\end{table}

These tests are pointwise in position and time. Spatially varying coefficients also contribute derivatives to the geodesic equations. If environmental inputs are updated between intervals, the metric and its forward distance are treated as time-independent within each interval; this piecewise-static procedure does not provide an exact solution for continuously varying environmental conditions.

\bibliography{bib}

@book{arnol2013mathematical,
  title={Mathematical methods of classical mechanics},
  author={Arnold, V. I.},
  volume={60},
  year={New York, Springer, 2013},
 publisher={Springer Science \& Business Media}
}

@phdthesis{dehkordi2018finsler,
  title={{F}insler Transnormal functions and singular foliations of codimension 1},
  author={Dehkordi, H. R.},
  year={IME, University of Sao paulo, Brazil, 2018},
  school={}
}

@article{alexandrino2019finsler,
  title={On Finsler transnormal functions},
  author={Alexandrino, Marcos M and Alves, Benigno O and Dehkordi, Hengameh R},
  journal={Differential Geometry and its Applications},
  volume={65},
  pages={93--107},
  year={2019},
  publisher={Elsevier}
}

@article{dehkordi2019huygens,
  title={Huygens’ envelope principle in {F}insler spaces and analogue gravity},
  author={Dehkordi, H. R. and Saa, A.},
  journal={Classical Quant. Grav.},
  volume={36}, 
  number={8},
  pages={085008},
  year={2019},
  publisher={IOP Publishing}
}

@article{markvorsen2016finsler,
  title={A Finsler geodesic spray paradigm for wildfire spread modelling},
  author={Markvorsen, S.},
  journal={Nonlinear Anal. Real World Appl.},
  volume={28}, 
  pages={208--228},
  year={2016},
  publisher={Elsevier}
}

@article{richards1990elliptical,
  title={An elliptical growth model of forest fire fronts and its numerical solution},
  author={Richards, Gwynfor D},
  journal={International Journal for Numerical Methods in Engineering},
  volume={30},
  number={6},
  pages={1163--1179},
  year={1990},
  publisher={Wiley Online Library}
}

@article{javaloyes2021applications,
  title={Applications of cone structures to the anisotropic rheonomic Huygens’ principle},
  author={Javaloyes, Miguel {\'A}ngel and Pend{\'a}s-Recondo, Enrique and S{\'a}nchez, Miguel},
  journal={Nonlinear Analysis},
  volume={209},
  pages={112337},
  year={2021},
  publisher={Elsevier}
}

@techreport{andrews2018rothermel,
  title={The {R}othermel surface fire spread model and associated developments: A comprehensive explanation},
  author={Andrews, Patricia L.},
  year={2018},
  number={RMRS-GTR-371},
  type={General Technical Report},
  institution={U.S. Department of Agriculture, Forest Service, Rocky Mountain Research Station},
  address={Fort Collins, CO}
}

@article{dehkordi2022applications,
  title={Applications of {R}anders geodesics for wildfire spread modelling},
  author={Dehkordi, Hengameh R.},
  journal={Applied Mathematical Modelling},
  volume={106},
  pages={45--59},
  year={2022},
  doi={10.1016/j.apm.2022.01.021},
  publisher={Elsevier}
}

@article{sharples2009overview,
  title={An overview of mountain meteorological effects relevant to fire behaviour and bushfire risk},
  author={Sharples, Jason J.},
  journal={International Journal of Wildland Fire},
  volume={18},
  number={7},
  pages={737--754},
  year={2009},
  doi={10.1071/WF08041}
}

@article{markvorsen2017geodesic,
  title={Geodesic sprays and frozen metrics in {R}heonomic {L}agrange manifolds},
  author={Markvorsen, S.},
  journal={arXiv preprint arXiv:1708.07350},
  year={2017}
}

@article{angel2023general,
  title={A general model for wildfire propagation with wind and slope},
  author={Javaloyes, Miguel {\'A}ngel and Pend{\'a}s-Recondo, Enrique and S{\'a}nchez, Miguel},
  journal={SIAM Journal on Applied Algebra and Geometry},
  volume={7},
  number={2},
  pages={414--439},
  year={2023},
  publisher={SIAM}
}

@article{dehkordi2023applications,
  title={Applications of {R}anders metric to track the paths through which the fire surrounds the power transmission lines},
  author={Dehkordi, Hengameh R},
  journal={Proceeding Series of the Brazilian Society of Computational and Applied Mathematics},
  volume={10},
  number={1},
  pages={2--7},
  year={2023}
}

@book{speer2022wildland,
  title={Wildland fire dynamics: fire effects and behavior from a fluid dynamics perspective},
  author={Speer, Kevin and Goodrick, Scott},
  year={2022},
  publisher={Cambridge University Press}
}

@book{rothermel1972mathematical,
  title={A mathematical model for predicting fire spread in wildland fuels},
  author={Rothermel, Richard C},
  volume={115},
  year={1972},
  publisher={Intermountain Forest \& Range Experiment Station, Forest Service, US}
}

@article{farsite2004fire,
  title={Fire area simulator-model development and evaluation},
  author={Finney, Mark A.},
  journal={Research Paper RMRS-RP-4 Revised. Ogden, UT: US Department of Agriculture, Forest Service, Rocky Mountain Research Station},
  year={2004}
}

@article{Sullivan2009,
author = {Sullivan, A L},
journal = {IJWF},
number = {4},
pages = {387--403},
publisher = {CSIRO},
title = {{Wildland surface fire spread modelling, 1990–2007. 3: Simulation and mathematical analogue models}},
volume = {18},
year = {2009}
}

@incollection{bressan2013,
  title={Dynamic blocking problems for a model of fire propagation},
  author={Bressan, A.},
  booktitle={Advances in Applied Mathematics, Modeling, and Computational Science},
  pages={11--40},
  year={2013},
  publisher={Springer}
}

@article{dehkordi2023wildfire,
  title={On the Wildfire Propagation in Flat Lands},
  author={Dehkordi, HR},
  journal={Trends in Computational and Applied Mathematics},
  volume={24},
  pages={83--90},
  year={2023},
  publisher={SciELO Brasil}
}

@article{dai1999diurnal,
  title={Diurnal and semidiurnal variations in global surface wind and divergence fields},
  author={Dai, Aiguo and Deser, Clara},
  journal={Journal of Geophysical Research: Atmospheres},
  volume={104},
  number={D24},
  pages={31109--31125},
  year={1999},
  publisher={Wiley Online Library}
}

@article{dehkordi2026satellite,
  author  = {Dehkordi, Hengameh R. and Kumar, Amit},
  title   = {Satellite-assisted wildfire propagation under
             environmental uncertainty: Huygens-based front evolution},
  journal = {Scientific Reports},
  year    = {2026},
  doi     = {10.1038/s41598-026-80574-2},
  url     = {https://doi.org/10.1038/s41598-026-80574-2}
}

@incollection{dehkordi2026determining,
  title={Determining a reliable model for wildfire propagation},
  author={Dehkordi, Hengameh R},
  booktitle={Reliability Assessment of Flow Networks and Industrial Systems},
  pages={205--216},
  year={2026},
  publisher={Elsevier}
}

@article{dehkordi2021mathematical,
  title={Mathematical modeling of wildfire propagation in an agricultural land},
  author={Dehkordi, Hengameh R},
  journal={Proceeding Series of the Brazilian Society of Computational and Applied Mathematics},
  volume={8},
  number={1},
  year={2021}
}
\bibliographystyle{unsrt}

\end{document}